\documentclass[12pt]{amsart}
\usepackage{amsmath,amsthm,upref,amssymb}
\usepackage[dvips]{graphicx}
\usepackage{psfrag}
\usepackage{wrapfig}
\usepackage[all]{xy}

\makeatletter
\renewcommand\paragraph%
  {\@startsection{paragraph}{4}{\z@}%
  {1.25ex \@plus1ex \@minus.2ex}{-1em}%
  {\normalfont\normalsize\bfseries}}
\makeatother

\newcommand*\Z{\mathbf Z}

\newcommand*\R{\mathbf R}

\newcommand*\T{\mathbf T}
\newcommand*\HH{\mathbf H}

\newcommand*\cH{{\mathcal{H}}}

\newcommand*\cU{{\mathcal{U}}}

\newcommand*\cM{{\mathcal{M}}}

\DeclareMathOperator{\lev}{lev}
\DeclareMathOperator{\lcm}{lcm}
\DeclareMathOperator\SL{SL}
\DeclareMathOperator\Lsurf{L}

\newcommand*{\lquotient}[2]{\left.\raisebox{-0.2ex}{$#2$}%
  \backslash\raisebox{0.2ex}{$#1$}\right.}

\newcommand\mat[4]{
  \bigl( \begin{smallmatrix} #1&#3\\ #2&#4 
  \end{smallmatrix}\bigr)
  }
\newcommand\hp[1]{
  \mat{1}{0}{#1}{1}%
  }
\newcommand\vp[1]{
  \mat{1}{#1}{0}{1}%
  }

\newtheorem{Theorem}{Theorem}
\newtheorem{ExternalTheorem}{Theorem}
\newtheorem{Conjecture}{Conjecture}
\newtheorem{Proposition}{Proposition}[section]
\newtheorem{Lemma}[Proposition]{Lemma}
\newtheorem{Corollary}[Proposition]{Corollary}
\theoremstyle{remark}
\newtheorem*{Remark}{Remark}
\newtheorem*{Thanks}{Acknowledgements}
\begin{document}

\title
{Noncongruence subgroups in $\cH(2)$}

%
\author{Pascal Hubert}
\address{%
IML, UMR CNRS 6206, Universit\'e de la M\'editerran\'ee,
Campus de Luminy, case 907, 13288 Marseille cedex 9, France%
}
\email{hubert@iml.univ-mrs.fr}
\author{Samuel Leli\`evre}
\address{%
\textsc{Irmar}, UMR CNRS 6625, Universit\'e de Rennes~1,
Campus Beaulieu, 35042 Rennes cedex, France;%
\newline
\indent
I3M, UMR CNRS 5149, Universit\'e Montpellier~2, case 51,
Place Eug\`ene Bataillon, 34095 Montpellier cedex 5, France;%
\newline
\indent
IML, UMR CNRS 6206, Universit\'e de la M\'editerran\'ee,
Campus de Luminy, case 907, 13288 Marseille cedex 9, France.%
}
\email{samuel.lelievre@polytechnique.org}
\urladdr{http://carva.org/samuel.lelievre/}

\date{29 May 2004}

\begin{abstract}
{We study the congruence problem for subgroups of the modular group 
that appear as Veech groups of square-tiled surfaces in the minimal 
stratum of abelian differentials of genus two.}
\end{abstract}

\maketitle
\thispagestyle{empty}

{\small Keywords: {congruence problem, Veech group,  
square-tiled surfaces}


\footnotesize
\setcounter{tocdepth}{1}
\tableofcontents
\normalsize

%

\section{Introduction}

Let $\omega$ be a holomorphic $1$-form on a compact Riemann surface
$X$.  If there exists a branched covering $f:X\to \T^2=\R^2/\Z^2$,
ramified only over the origin of $\T^2$, such that $f^* (dz)=\omega$,
the flat surface $(X,|\omega|)$ is tiled by squares whose vertices
project to the origin of the torus, and $(X,\omega)$ is called a
\textbf{square-tiled (translation) surface}.

In each genus $g$, square-tiled surfaces are the integer points of the
moduli space $\cH_g = \Omega\cM_g$ of holomorphic $1$-forms on Riemann
surfaces of genus $g$.  This space is stratified by the combinatorial
type of zeros, and each stratum is a complex orbifold endowed with an
action of $\SL(2,\R)$.  Orbits for this action are called
Teichm\"uller discs.

The main problem in dynamics in Teichm\"uller spaces is to understand
this $\SL(2,\R)$-action, and to obtain Ratner-like classification
results for its orbit closures and its invariant closed submanifolds.

The first step is to determine as many invariant closed submanifolds 
as possible.
The simplest of them are closed orbits.  These are the orbits of
translation surfaces with finite-covolume stabilisers, called Veech
surfaces because of Veech's pioneering work \cite{Ve}.  These
Teichm\"uller discs project to geodesically embedded curves, called
Teichm\"uller curves, in the moduli space $\cM_g$ of complex curves of
genus $g$.  These curves are uniformised by the stabiliser of the 
corresponding $\SL(2,\R)$-orbit.

Square-tiled surfaces are Veech surfaces.  They already appeared in
Thurston's work on the classification of surface diffeomorphisms, see
\cite[expos\'e~13]{FLP}.  Nevertheless up to recently their
Teichm\"uller discs have been little discussed, due to the difficulty
of proving precise statements about them.  The only classical result
is Gutkin and Judge's theorem \cite{GuJu} which states that the
corresponding stabilisers are arithmetic (commensurable to
$\SL(2,\Z)$).  Very recently the Teichm\"uller discs of square-tiled
surfaces were studied into more detail, see \cite{HL}, \cite{Mc4}, 
\cite{Moe}, \cite{Schmi}.

A square-tiled surface $(X,\omega)$ is called \textbf{primitive} if
the lattice of relative periods of $\omega$ is $\Z^2$ (in other words
the covering $(X,\omega)\to (\T^2,dz)$ does not factor through a
bigger torus).  In this case, the stabiliser, denoted by
$\SL(X,\omega)$, is a (finite-index) subgroup of $\SL(2,\Z)$.

In order to give the most accurate description of Teichm\"uller discs
of square-tiled surfaces, we investigate these subgroups.  In the
theory of subgroups of $\SL(2,\Z)$, a natural and important question
is the congruence problem.  This question is the central object of
this paper: we give a negative answer in the stratum
$\cH(2)=\Omega\cM_2(2)$ of $1$-forms on genus $2$ surfaces having one
double zero.


\paragraph{Recent results about square-tiled surfaces in $\cH(2)$.}

The discrete orbit $\SL(2,\Z){\cdot}(X,\omega)$ of a primitive
square-tiled surface $(X,\omega)$ consists of all the primitive
square-tiled surfaces in its Teichm\"uller disc
$\SL(2,\R){\cdot}(X,\omega)$; indeed, $\SL(2,\Z)$ acts on primitive
square-tiled surfaces, preserving the number of squares.
Understanding the Teichm\"uller discs or the discrete orbits of
primitive square-tiled surfaces is therefore equivalent.  We will use
the following result about the discrete orbits of primitive
square-tiled surfaces in $\cH(2)$.

\begin{ExternalTheorem}
\label{thm:orbits}
%
Primitive $n$-square-tiled surfaces in the stratum $\cH(2)$ form: one
orbit $A_3$ if $n = 3$; two orbits $A_n$ and $B_n$ if $n$ is odd
$\geqslant 5$; one orbit $C_n$ if $n$ is even.
\end{ExternalTheorem}

This was shown for prime $n$ in \cite{HL}, and conjectured for 
arbitrary $n$; the conjecture was proved in full generality in 
\cite{Mc4}.

Let $\Gamma_{A_n}$, $\Gamma_{B_n}$ and $\Gamma_{C_n}$ denote the
stabilisers of these orbits.  


\begin{Remark}
The indices of the groups $\Gamma_{A_n}$, $\Gamma_{B_n}$,
$\Gamma_{C_n}$ in $\SL(2,\Z)$ are the cardinalities $a_n$, $b_n$,
$c_n$ of the discrete orbits $A_n$, $B_n$, $C_n$.
\end{Remark}

Eskin--Masur--Schmoll \cite{EsMaSc} give a formula for the number of
primitive $n$-square-tiled surfaces in $\cH(2)$:

\begin{ExternalTheorem}
\label{thm:global:primitive:count}
The number of primitive $n$-square-tiled surfaces in $\cH(2)$ is
$\frac{3}{8}(n-2)n ^2\prod_{p\mid n}(1-\frac{1}{p^2})$.
\end{ExternalTheorem}

\begin{Remark}
Throughout this paper, the letter $p$ always denotes prime
numbers; in particular, $\prod_{p\mid n}$ is the product over prime
divisors of $n$.
\end{Remark}

This formula gives $c_n$ (and $a_3$) when there is one orbit and $a_n
+ b_n$ when there are two.  We conjectured in \cite{HL}:

\begin{Conjecture}
\label{conj:index}
For odd $n\geqslant 5$, $a_n$ and $b_n$ are given by:

\centerline{$a_n = \frac{3}{16}(n-1)n ^2\prod_{p\mid
n}(1-\frac{1}{p^2}), \quad b_n = \frac{3}{16}(n-3)n ^2\prod_{p\mid
n}(1-\frac{1}{p^2})$.}
\end{Conjecture}

\paragraph{Statement of results.}

In this paper, we show:

\begin{Theorem}
\label{thm:noncongruence}
For all even $n\geqslant 4$, $\Gamma_{C_n}$ is a noncongruence
subgroup.  For all odd $n\geqslant 5$ satisfying
Conjecture \ref{conj:index}, $\Gamma_{A_n}$ and $\Gamma_{B_n}$
are noncongruence subgroups.
\end{Theorem}

\begin{Remark}
Conjecture \ref{conj:index} is proved up to $n = 10000$ by an explicit 
combinatorial computer calculation.
\end{Remark}

\begin{Corollary} Under Conjecture \ref{conj:index}, the only
primitive square-tiled surfaces in $\cH(2)$ whose stabiliser is a
congruence subgroup are those tiled with $3$ squares.
\end{Corollary}

\begin{Corollary}
Under Conjecture \ref{conj:index}, of all the Teichm\"uller curves
embedded in $\cM_2$ that come from orbits in $\cH(2)$, only one is
uniformised by a congruence subgroup of $\SL(2,\Z)$.
\end{Corollary}

\begin{Remark}
For $n = 3$, $\Gamma_{A_3}$ is the level $2$ congruence subgroup 
$\Theta$ generated by $\mat{0}{1}{-1}{0}$ and $\hp{2}$, named 
after its link to the Jacobi Theta function.
\end{Remark}

\paragraph{Link with the Hurwitz problem.}
An essential ingredient in our proof of Theorem
\ref{thm:noncongruence} is the knowledge of the indices in $\SL(2,\Z)$
of $\Gamma_{A_n}$, $\Gamma_{B_n}$ and $\Gamma_{C_n}$ (given by Theorem
\ref{thm:global:primitive:count} and Conjecture \ref{conj:index}).

Since these indices are the cardinalities of the discrete orbits
$A_n$, $B_n$ and $C_n$, finding these numbers is a variant of
Hurwitz's problem, which consists in counting the number of branched
covers of a fixed combinatorial type (number and multiplicity of
ramification points) and fixed degree of a Riemann surface $S$.  A
very detailed survey of this subject can be found in the introduction
of Zvonkine's thesis \cite{Zv}.

When $S$ is the torus $\T^2 = \R^2/\Z^2$ (or more generally an
elliptic curve), it can be endowed with the $1$-form $dz$.  Hurwitz's
problem amounts to counting the number of coverings (with fixed
combinatorial type) $f:(X,\omega)\to(\T^2,dz)$ where $\omega=f^*(dz)$.
For a fixed combinatorial type $c$, denote by $h_{n,c}$ the number of
such coverings, weighted by the inverse of their number of
automorphisms.

We have the following fundamental theorem:

\begin{ExternalTheorem}
For any combinatorial type, the generating series 
$F_c(z)= \sum_{h=1}^\infty h_{n,c}\,q^n$, where $q = e^{2i\pi z}$, 
is a quasi-modular form of maximal weight $6g-6$.
\end{ExternalTheorem}

This theorem was first proved in the case of simple ramifications by
Dijkgraaf \cite{Di} and Kaneko--Zagier \cite{KaZa}; the general proof
relies on results of Bloch--Okounkov \cite{BlOk}, see \cite{EsOk}.

The quasi-modular form is explicitated by Kani \cite{Ka} and by 
Eskin--Masur--Schmoll \cite{EsMaSc} in particular cases. Some 
generalisations are proved by Eskin--Okounkov--Pandharipande 
\cite{EsOkPa}.

Note also that the asymptotics of the countings of square-tiled
surfaces of bounded area serve to compute the volumes of strata (see
\cite{Zo}, \cite{EsOk}).

\begin{Thanks}
We thank Gabriela Schmith\"usen for the inspiration and useful 
discussions, and Giovanni Forni for encouraging this research.
\end{Thanks}

\section{Background}

\subsection{Square-tiled surfaces, action of $\SL(2,\Z)$, cusps}

We recall here some tools used in \cite{HL}, to which we refer for 
more detail.

The modular group $\Gamma(1) = \SL(2,\Z)$ acts on primitive
square-tiled surfaces, preserving the number of squares tiles.
Indeed, the property of having $\Z^2$ as lattice of relative periods
is $\SL(2,\Z)$-invariant.

Given a primitive square-tiled surface $(X,\omega)$, its stabiliser
$\SL(X,\omega)$ is a finite-index subgroup of $\SL(2,\Z)$,
therefore the curve $\lquotient{\HH}{\SL(X,\omega)}$ is a branched
cover of the modular curve $\lquotient{\HH}{\SL(2,\Z)}$, and the 
degree of the cover is the index of $\SL(X,\omega)$ in $\SL(2,\Z)$.

The modular group is generated by any two
matrices among $\hp{1}$, $\vp{1}$ and $\mat{0}{1}{-1}{0}$.  Denote by
$\cU$ the subgroup generated by $\hp{1}$.

\paragraph{Cusps.}
The cusps of $\lquotient{\HH}{\SL(X,\omega)}$ are classified 
combinatorially by the following lemma.

\begin{Lemma}[Zorich]
Let $(X,\omega)$ be a primitive square-tiled surface.  There is a
$1$-$1$ correspondence between the set of cusps of
$\lquotient{\HH}{\SL(X,\omega)}$ and the $\cU$-orbits of 
$\SL(2,\Z){\cdot}(X,\omega)$.
\end{Lemma}

Any square-tiled surface decomposes into horizontal cylinders, which
are also square-tiled, and bounded by unions of saddle connections of
integer lengths.  This provides a way to give coordinates for
square-tiled surfaces in each stratum (see below for the stratum 
$\cH(2)$).

The action of the generators $\hp{1}$ and $\mat{0}{1}{-1}{0}$ of
$\SL(2,\Z)$ is easily seen in these coordinates:
$\mat{0}{1}{-1}{0}$ exchanges the horizontal and vertical
directions; $\hp{1}$ only changes the twists.

The width of a cusp is given by the cardinality of the corresponding
$\cU$-orbit. If the horizontal cusp has width $\ell$, the primitive 
parabolic in the horizontal direction is $\hp{\ell}$. Considering how 
the cylinders behave under the action of $\cU$, we get the following 
lemma.

\begin{Lemma}
\label{lem:cusp:width}
If a primitive square-tiled surface decomposes into horizontal
cylinders $c_i$ of height $h_i$ and width $w_i$, then its (horizontal)
cusp width equals the least common multiple of the $\frac{w_i}{h_i
\wedge w_i}$, possibly divided by some factor.
\end{Lemma}

\paragraph{Notation.}
Here, and in the sequel, $a\wedge b$ denotes the greatest common 
divisor of two integers $a$ and $b$.

\medskip

The following example illustrates the case of division by a factor.

\begin{wrapfigure}{r}{0pt}
\psfrag{a}{\Small $a$}
\psfrag{b}{\Small $b$}
\psfrag{c}{\Small $c$}
\psfrag{d}{\Small $d$}
\includegraphics{smallcw.eps}
\end{wrapfigure}

This surface is in $\cH(1,1)$ and has a nontrivial translation by the
vector $(2,0)$; though it is made of one cylinder of height $1$ and
width $4$, its cusp width is only $2$.

In the stratum $\cH(2)$ on which we will focus from now on, this 
situation does not occur.

\subsection{Square-tiled surfaces in $\cH(2)$}
The stratum $\cH(2)$ has recently received much attention
(\cite{EsMaSc}, \cite{Ca}, \cite{Mc1,Mc3,Mc4}, \cite{HL}).

Square-tiled surfaces in $\cH(2)$ are of two types \cite{Zo}, the
one-cylinder ones and the two-cylinder ones.  The corresponding
coordinates are: for one-cylinder surfaces, one height, three lengths
of saddle connections and one twist parameter; for two-cylinder
surfaces, one height, width and twist for each cylinder.

Theorem \ref{thm:orbits} says that for each odd $n \geqslant 5$,
primitive $n$-square-tiled surfaces are in two orbits $A_n$ and $B_n$.
These orbits are distinguished by a simple invariant, the number
of integer Weierstrass points (i.e.\ Weierstrass points located at
vertices of the square tiles).  A surface is in $A_n$ if it has one
integer Weierstrass point, in $B_n$ if it has three.

The coordinates for square-tiled surfaces in $\cH(2)$ were used in
\cite{Zo}, in \cite{EsMaSc} and in \cite{HL} where the position of 
Weierstrass points was also discussed and the invariant introduced.
This invariant was independently expressed in terms of divisors by
Kani \cite{Ka}.  McMullen \cite{Mc4} expressed it as the parity of a
spin structure.

\begin{wrapfigure}{r}{0pt}
\psfrag{h1}{\Small $h_1$}
\psfrag{h2}{\Small $h_2$}
\psfrag{w1}{\Small $w_1$}
\psfrag{w2}{\Small $w_2$}
\psfrag{t1}{\Small $t_1$}
\psfrag{t2}{\Small $t_2$}
\includegraphics{wpts2.eps}
\end{wrapfigure}

\paragraph{Notation}Denote by $S(h_1, h_2, w_1, w_2, t_1, 
t_2)$ the two-cylinder surface with cylinders $c_i$ of height $h_i$, 
width $w_i$ and twist $t_i$, with $w_1 < w_2$.

The figure shows a fundamental polygon for $S(2,3,3,8,2,1)$; the
surface is obtained from this polygon by identifying pairs of parallel
sides of same lengths.  We indicate the double zero by black dots and
the other Weierstrass points by circles.  The same conventions hold
for all pictures in this paper.

Let us give some examples of square-tiled surfaces in $\cH(2)$.

First, some one-cylinder surfaces of particular interest.

\begin{Lemma}
For each $n\geqslant 4$, there is a primitive $n$-square-tiled 
surface which is one-cylinder both horizontally and vertically.
\end{Lemma}

\begin{wrapfigure}{r}{0pt}
\psfrag{1}{\Small $1$}
\psfrag{n}{\Small $n-3$}
\psfrag{2}{\Small $2$}
\includegraphics{hvonecyl}
\end{wrapfigure}

The one-cylinder surface with saddle connections of 
lengths $1$, $n-3$, $2$ on the top and $2$, $n-3$, $1$ on the bottom 
has this property.

\begin{Corollary}
The stabiliser of this surface contains $\hp{n}$ and $\vp{n}$.
\end{Corollary}

Indeed, one-cylinder cusps have width $n$.

\begin{Remark}
When $n$ is odd, the surface described above is in orbit $B_n$.
\end{Remark}

Some two-cylinder surfaces also deserve special attention.

\begin{wrapfigure}{r}{0pt}
\psfrag{1}{\Small $1$}
\psfrag{a}{\Small $a$}
\psfrag{b}{\Small $b$}
\includegraphics{Lab.eps}
\end{wrapfigure}

\paragraph{Notation}
For $a$ and $b$ $\geqslant 2$, denote by $\Lsurf(a,b)$ the surface
$S(a-1,1,1,b,0,0)$.  This surface is a primitive square-tiled 
surface tiled by $n = a + b - 1$ squares.
This surface has cusp width $b$ and vertically $a$.

When $n$ is odd, this surface is in $A_n$ if $a$ and $b$ are even, 
in $B_n$ if $a$ and $b$ are odd.

\subsection{Congruence subgroups; level of a subgroup}

The material in this section is classical, and can be found in \cite{Ra}. 

For any integer $m>1$, consider the natural projection $\SL(2,\Z) \to
\SL(2,\Z/m\Z)$.  This projection is a group homomorphism.  Its kernel
is called the \textbf{principal congruence subgroup of level $m$}, and denoted
by $\Gamma(m)$.  It consists in all matrices congruent to
$\mat{1}{0}{0}{1}$ modulo $m$. 
This is consistent with the notation
$\Gamma(1)$ for $\SL(2,\Z)$.

\begin{Lemma}
For any $m$, $[\Gamma(1):\Gamma(m)]=m^3\prod_{p\mid
m}(1-\frac{1}{p^2})$.
\end{Lemma}

\begin{Corollary}
\label{cor:index:congruence:subgroups}
%
If $m\wedge m' = 1$, then $[\Gamma(m):\Gamma(mm')] =
[\Gamma(1):\Gamma(m')]$.
\end{Corollary}

Any group $\Gamma$ containing some $\Gamma(m)$ is called a \textbf{congruence
subgroup}, and its \textbf{level} is defined to be the least $m$ such that 
$\Gamma(m) \subset \Gamma$ (i.e.\ the level of the largest principal 
congruence subgroup it contains).

\begin{Remark}
A principal congruence subgroup is a normal subgroup of $\Gamma(1)$.
Hence being a congruence subgroup is invariant by conjugation in
$\SL(2,\Z)$; the level is also invariant.
\end{Remark}

There is a more general notion of level, due to Wohlfahrt \cite{Wo}.
The \textbf{level} of a finite-index subgroup of $\SL(2,\Z)$ is the least
common multiple of its cusp widths.  Wohlfahrt proved that for
congruence subgroups, it coincides with the previous definition, and
that:

\begin{Lemma}[Wohlfahrt \cite{Wo}]
\label{lem:cg:iff:lev}
A finite-index subgroup of level $\ell$ is a congruence subgroup if
and only if it contains the principal congruence subgroup of level
$\ell$.
\end{Lemma}

\subsection{Quasi-modular forms}

As said in the introduction, the generating function for the weighted
countings of surfaces tiled by $n$ squares is a quasi-modular form.

The numbers $h_{n,c}$ of surfaces tiled by $n$ squares in a given
stratum, and the numbers $h_{n,c}^\text{P}$ of primitive ones, 
are related by 

\noindent\centerline{$h_{n,c} = \sum_{d\mid n} \sigma(n/d) h_{d,c}^\text{P}$,}

\noindent where $\sigma(k) = \sum_{d\mid k} d$ is the sum of divisors
of $k$.  This is because the number of tori tiled by $n$ squares is
$\sigma(n)$.

In addition, we note that in $\cH(2)$, the coverings have no 
automorphisms, hence the weighted and unweighted countings are the same.

\begin{Conjecture}
\label{conj:quasimod}
In $\cH(2)$, the countings for odd $n$ according to the invariant are
generated by a quasi-modular form.
\end{Conjecture}

Theorem \ref{thm:global:primitive:count} is mentioned in
\cite{EsMaSc} as a consequence of the quasi-modularity.
Likewise, Conjecture \ref{conj:index} would follow from Conjecture 
\ref{conj:quasimod}.


\newpage
\section{Strategy for the proof of Theorem \ref{thm:noncongruence}}

We build on the proof by Schmith\"usen \cite{Schmi} that the
stabiliser of a $4$-square-tiled surface in $\cH(2)$ is a
noncongruence subgroup, based on an idea of Stefan K\"uhnlein.

\subsection{Sufficient conditions for noncongruence}
\label{sec:conditions}

%
\setlength{\intextsep}{-3pt}
\begin{wrapfigure}{r}{0pt}
\footnotesize
%
$\xymatrix@C=-4ex@R=2.5ex@M=5pt{
& \Gamma(1)=\SL(2,\Z) & \\
\Gamma\ar@{^{(}->}[ur]_{\textstyle d}& & 
\Gamma(m)\ar@{_{(}->}[ul]\\
& \Gamma \cap \Gamma(m)\ar@{_{(}->}[ul]\ar@{^{(}->}[ur]^{\textstyle d'} & \\
& & \Gamma(\ell)\ar@{^{(}->}[uu]^{\textstyle \delta}}$
\end{wrapfigure}

\normalsize
\setlength{\intextsep}{0pt}

Let $\Gamma$ be a subgroup of $\Gamma(1)$ of finite
index $d$ and level $\ell$.  For any divisor $m$ of $\ell$, consider
the finite-index inclusions represented on the figure.

\paragraph{Two remarks.} First, if $\Gamma$ projects surjectively to
$\SL(2,\Z/m\Z)$, one can conclude by observing the two
exact sequences below that $d' = d$, where $d' =
[\Gamma(m):\Gamma\cap\Gamma(m)]$ and $d = [\Gamma(1):\Gamma]$.

\centerline{%
\footnotesize
$\xymatrix@R=3ex@M=5pt{1\ar@{=}[d]\ar[r] & \Gamma(m)\ar[r] &
\Gamma(1)\ar[r] & \SL(2,\Z/m\Z)\ar@{=}[d]\ar[r] & 
1\ar@{=}[d]\\
1\ar[r] & 
\Gamma\cap\Gamma(m)\ar@{^{(}->}[u]_{\textstyle{d'}}\ar[r] 
& \Gamma\ar@{^{(}->}[u]_{\textstyle{d}}\ar[r] & 
\SL(2,\Z/m\Z)\ar[r] & 1}$%
}%
\normalsize

Second, if $\Gamma$ is a congruence subgroup, and hence by Lemma
\ref{lem:cg:iff:lev} contains $\Gamma(\ell)$, then
$\Gamma(\ell)$ is contained in $\Gamma\cap\Gamma(m)$ and the indices satisfy
$[\Gamma(m):\Gamma(\ell)] = [\Gamma(m):\Gamma\cap\Gamma(m)] {\cdot}
[\Gamma\cap\Gamma(m):\Gamma(\ell)]$, which implies $d'\mid \delta$.

Combining these two remarks, we get the following sufficient condition
for noncongruence, which was used by Schmith\"usen \cite{Schmi}.

\begin{Proposition}[K\"uhnlein]
If $\Gamma$ is a subgroup of $\Gamma(1)$ of finite index $d$ and 
level $\ell$ and there exists a divisor $m$ of $\ell$ for which\\
\textbullet\ $\Gamma$ projects surjectively to $\SL(2,\Z/m\Z)$, and\\
\textbullet\ the index $\delta = [\Gamma(m):\Gamma(\ell)]$ is not
a multiple of $d$,\\
then $\Gamma$ is not a congruence subgroup.
\end{Proposition}

\begin{Remark}
Suppose $\Gamma$ contains two matrices $\hp{k}$ and
$\vp{k'}$.  If $m$ is an integer relatively prime to both
$k$ and $k'$, then $k$ and $k'$ are invertible modulo $m$ so some
powers of $\hp{k}$ and $\vp{k'}$ project to
$\hp{1}$ and $\vp{1}$ in $\SL(2,\Z/m\Z)$, hence 
the projection $\GammaÊ\to \SL(2,\Z/m\Z)$ is surjective.
\end{Remark}

This extra remark yields the following sufficient condition for 
noncongruence.

\begin{Proposition}
\label{prn:crit:ncg}
If a subgroup $\Gamma\subset\Gamma(1)$ of finite index $d$ contains
two matrices $\hp{k}$ and $\vp{k'}$ and if its
level $\ell$ has a divisor $m$ relatively prime to both $k$ and $k'$,
such that the index $\delta = [\Gamma(m):\Gamma(\ell)]$ is not a
multiple of $d$, then $\Gamma$ is not a congruence subgroup.
\end{Proposition}

\subsection{Strategy}

Consider an orbit $A_n$, $B_n$ or $C_n$, with $n$ as in
Theorem~\ref{thm:noncongruence}.  Its stabiliser $\Gamma_{A_n}$,
$\Gamma_{B_n}$ or $\Gamma_{C_n}$ is defined only up to conjugation in
$\SL(2,\Z)$; the representatives of the conjugacy class are the
stabilisers of the (square-tiled) surfaces in the orbit.  The index and level are
preserved by conjugation in $\SL(2,\Z)$.


\paragraph{Choice.}
Let $S$ be a (square-tiled) surface in an orbit $A_n$, $B_n$ or $C_n$, and $\Gamma$
be its stabiliser.  

\paragraph{Notation.}
Denote by $d$ the index of $\Gamma$ and by $\ell$
its level.  Consider the prime factor decompositions $n = \prod p^\nu$
and $\ell = \prod p^\lambda$, where $\nu$ and $\lambda$ can denote a
different integer for each prime $p$.

\paragraph{Choice.}
Choose some $\hp{k}$ and $\vp{k'}$ in $\Gamma$, for instance $k$ and
$k'$ could be taken to be the horizontal and vertical cusp widths of
$S$.

\paragraph{Notation}
Following \cite{Mc4}, if $a$ and $b$ are two integers, denote by 
$a /\!/ b$ the greatest divisor of $a$ that is prime to $b$.
If $a = \prod p^\alpha$ is the prime factor decomposition of $a$, we
have $a /\!/ b = \prod_{p\,\nmid\, b} p^\alpha = a / \prod_{p\mid b}
p^\alpha$.

\paragraph{Choice.}
Choose $m = \ell /\!/ kk' = \ell / \prod_{p\mid kk'}p^\lambda$.

\paragraph{Notation.}
Denote by $\delta$ the index of $\Gamma(\ell)$ in $\Gamma(m)$.

By construction $m$ is a divisor of $\ell$, relatively prime to both
$k$ and $k'$.  In view of applying Proposition \ref{prn:crit:ncg},
there remains only to check that $d$ does not divide $\delta$.
Since $m$ is also relatively prime to $\ell/m$, by
Corollary \ref{cor:index:congruence:subgroups}, $\delta = 
(\ell/m)^3 \prod_{p\mid \ell/m}(1-\frac{1}{p^2})$.

\begin{Remark}
If $a$ is an integer and $a = \prod p^\alpha$ is its prime
factor decomposition, one can rewrite $a^r\prod_{p\mid
a}(1-\frac{1}{p^2})$ as $\prod_{p\mid a}p^{r\alpha-2}(p^2-1)$.
Hence \\
\textbullet\ $\delta = \prod_{p\mid kk'}p^{3\lambda-2}(p^2-1)$, and \\
\textbullet\ $d = f(n)\prod_{p\mid n}p^{2\nu-2}(p^2-1)$, where $f(n)$ is one of
$\frac{3}{16}(n-1)$, $\frac{3}{16}(n-3)$, $\frac{3}{8}(n-2)$, 
according to whether orbit $A_n$, $B_n$ or $C_n$ is under 
consideration.
\end{Remark}

In order to complete the proof, there merely remains to describe how
to apply our strategy.

For this we need the levels of $\Gamma_{A_n}$, $\Gamma_{B_n}$ and
$\Gamma_{C_n}$; we give them in \S\,\ref{sec:levels}.

The last three sections then describe, in each orbit, good choices of
a surface $S$, values of $k$ and $k'$, and, keeping the notations
($d$, $\ell$, $\nu$, $\lambda$, $m$, $\delta$) introduced here (and
consistent with those in \S\,\ref{sec:conditions}), show that $d$ does
not divide $\delta$.

\section{The level of $\Gamma_{A_n}$, $\Gamma_{B_n}$ and
$\Gamma_{C_n}$}
\label{sec:levels}

As said above, the stabiliser of an $\SL(2,\Z)$-orbit of square-tiled 
surfaces is defined up to conjugacy in $\SL(2,\Z)$, but its level
is well-defined.

\begin{Proposition}
\label{prn:lev}
The groups $\Gamma_{A_n}$, $\Gamma_{B_n}$ and
$\Gamma_{C_n}$ have levels:

\centerline{$\lev \Gamma_{A_n} = d_n, \quad 
\lev \Gamma_{B_n} = d_n / 4, \quad
\lev \Gamma_{C_n} = d_n,$}

\noindent
where $d_n = \lcm(1,2,3,\ldots,n)$.
\end{Proposition}

\begin{Remark}
The prime factor decomposition of $d_n$ is $\prod_{p\leqslant
n}p^\tau$ where the exponents $\tau$ are the integers such that
$p^\tau \leqslant n < p^{\tau+1}$.
\end{Remark}

The remainder of this section is devoted to proving the proposition.

First recall that the level of $\GammaÊ\subset \SL(2,\Z)$ is defined 
as the least common multiple of the amplitudes of the cusps of $\Gamma$. 
When $\Gamma$ is the stabiliser of a primitive square-tiled surface $S$, 
its cusp widths are equivalently the horizontal cusp widths of the 
surfaces in the $\SL(2,\Z)$-orbit of $S$.

Recall also Lemma \ref{lem:cusp:width}.  If $S$ is tiled by $n$
squares, the widths of its cylinders are at most $n$, so the level
of $\Gamma$ divides $\lcm(1,2,3,\ldots,n)$.

Orbit $C_n$ (for even $n$) contains one-cylinder surfaces, which
have cusp width $n$, and, for all $a$ and $b$ such that $a + b = n +
1$ and $2 \leqslant a, b \leqslant n - 1$, two-cylinder surfaces
$\Lsurf(a,b)$, which have cusp width $b$. Hence, the level of $\Gamma_{C_n}$
is a multiple of, and therefore equals, $\lcm(1,2,3,\ldots,n)$.

\begin{wrapfigure}{r}{0pt}
\psfrag{a}{\Small $a$}
\psfrag{b}{\Small $b$}
\includegraphics{a9.eps}
\end{wrapfigure}

Orbit $A_n$ (for odd $n$) contains one-cylinder surfaces, which
have cusp width $n$, and, for all $a$ and $b$ such that $a + b = n$
and $1 \leqslant a < b \leqslant n - 1$, two-cylinder surfaces with
two cylinders of height $1$ and widths $a$ and $b$, which have cusp
width $\lcm(a,b)$. Hence, the level of $\Gamma_{A_n}$ is a multiple of, and
therefore equals, $\lcm(1,2,3,\ldots,n)$.

Orbit $B_n$ (for odd $n$) contains one-cylinder surfaces, which
have cusp width $n$, and, for all odd $a$ and $b$ such that $a + b = n
+ 1$ and $2 \leqslant a, b \leqslant n - 1$, two-cylinder surfaces
$\Lsurf(a,b)$, which have cusp width $b$. Hence, the level of $\Gamma_{B_n}$
is a multiple of $\lcm(1, 3, 5,\ldots, n)$.

Since $\lcm(1,2,3,\ldots,n)$ is a power of $2$ times $\lcm(1, 3,
5,\ldots, n)$, there remains only to determine the power of $2$ in the
level of $\Gamma_{B_n}$, i.e.\ the maximal power of $2$ that can arise
as a divisor of $\frac{w}{h \wedge w}$ for the height $h$ and the
width $w$ of a cylinder of a surface of $B_n$.

Let $\tau$ be the integer such that $2^\tau < n < 2^{\tau + 1}$.

There is at least one two-cylinder surface $S(h_1,2,w_1, 2^{\tau - 1},
t_1, t_2)$ with odd $t_2$ in $B_n$; such a surface satisfies
$\frac{w_2}{h_2 \wedge w_2} = 2^{\tau - 2}$.

Suppose a surface $S$ in $B_n$ has even cusp width $k = 2^t {\cdot}q$
with $q$ odd.  Then $S$ has two cylinders, and by the discussion in
\cite[\S \,5.1]{HL}, one cylinder has even width $w$ and even height
$h$, while the other has odd height $h'$ and odd width $w'$.  Since
$k = \lcm(\frac{w}{w\wedge h},\frac{w'}{w'\wedge h'})$ and 
$\frac{w'}{w'\wedge h'}$ is odd, $2^t$ divides $\frac{w}{h \wedge w}$.
But $h \geqslant 2$ and since $n = hw+h'w'$, $w < n / h \leqslant n /
2 < 2^\tau$, so $\frac{w}{h \wedge w} \leqslant w / 2 < 2^{\tau-1}$. 
Therefore $t \leqslant \tau - 2$.

\section{Noncongruence of $\Gamma_{C_n}$ for even $n\geqslant 4$}

\subsection{Case when $n-2$ is not a power of $2$}

We take $S$ to be the one-cylinder surface with saddle connections of
lengths $1$, $n-3$, $2$ on the top and $2$, $n-3$, $1$ on the bottom.

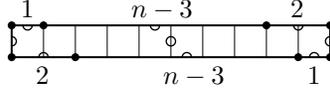
\begin{center}
\psfrag{1}{\Small $1$}
\psfrag{n}{\Small $n-3$}
\psfrag{2}{\Small $2$}
\includegraphics{even.eps}
\end{center}

As a one-cylinder surface it has cusp width $k = n$ and since its
vertical direction is also one-cylinder, its vertical cusp width $k'$
is also $n$, so $\Gamma$ contains $\hp{n}$ and $\vp{n}$.

Recall that $\Gamma$ has index $d = \frac{3}{8} (n-2) \prod_{p\mid n}
p^{2\nu-2} (p^2-1)$.

Choosing $m = \ell/\!/n = \ell / \prod_{p\mid n} p^\lambda$ leads to
$\delta = \prod_{p\mid n}p^{3\lambda-2}(p^2-1)$.

So $d$ divides $\delta$ if and only if $3(n-2)$ divides $2^{3} {\cdot}
\prod_{p\mid n} p^{3\lambda-2\nu}$.

Since $n \wedge (n - 2) = 2$, the assumption that $n - 2$ is not a
power of $2$ implies it has some (odd) prime factors that do not
divide $n$.

Hence $d$ does not divide $\delta$, so $\Gamma$ cannot be a congruence 
subgroup.

\subsection{Case when $n-2$ is a power of $2$}

The case $n = 4$ is known from \cite{Schmi}.  It can also be treated
as above, since the index of $\Gamma_{C_4}$ is $d = 9$ and, taking $S$
and $m$ as above, $\delta = 2^4{\cdot}3$.

\begin{wrapfigure}{r}{0pt}
\psfrag{n}{\Small $n-2$}
\includegraphics{even2.eps}
\end{wrapfigure}
From now on assume $n>4$.  

We take $S = S(1,1,1,n-2,1,0)$.

Note that this requires that $n-2 > 2$, which is why 
the case  $n=4$ was dealt with separately.

This surface has horizontal cusp width $n-2$
and vertical cusp width $4$, so the stabiliser $\Gamma$ contains
$\hp{n-2}$ and $\vp{4}$.

Recall that $\Gamma$ has index $d =
\frac{3}{8}(n-2)\prod_{p\mid n}p^{2\nu-2}(p^2-1)$.

Choosing $m = \ell/\!/2 = \ell/2^\lambda$ leads to $\delta =
2^{3\lambda-2}(2^2-1)=2^{3\lambda-2}{\cdot}3$.

Since $n$ is even, it has $p = 2$ as a prime factor, which gives $3$ 
as $p^2 - 1$, so $3^2$ divides $d$.

Hence $d$ does not divide $\delta$, so $\Gamma$ cannot be a congruence
subgroup.

\newpage
\section{Noncongruence of $\Gamma_{A_n}$ for odd $n\geqslant 5$}

\subsection{Case when $n - 1$ is a power of $2$}

Take $S= \Lsurf(2,n-1)$.

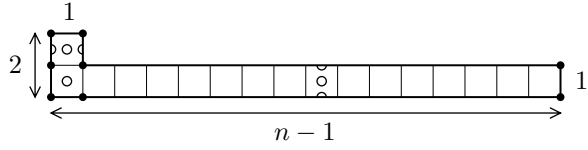
\begin{center}
\psfrag{1}{\Small $1$}
\psfrag{2}{\Small $2$}
\psfrag{n}{\Small $n-1$}
\includegraphics{af.eps}
\end{center}

Its cusp width is $n-1$ ($ = 2^\lambda$) and its vertical cusp width
is $2$, so its stabiliser $\Gamma$ contains $\hp{n-1}$ and $\vp{2}$.

Here $d = \frac{3}{16}(n-1)\prod_{p\mid n} p^{2\nu-2} (p^2-1)$.

The choice of $m = \ell /\!/ 2 = \ell/2^\lambda$ leads to
$\delta=2^{3\lambda-2} {\cdot} 3$.

If $n$ is a power of $3$, then $3^3$ divides $d$; otherwise $n$ has
some (odd) prime factor $p \neq 3$, for which $p^2 - 1 = (p-1)(p+1)$
is a multiple of $3$, so that $3^2$ divides $d$.  Therefore $d$ does
not divide $\delta$ and $\Gamma$ is not a congruence subgroup.

\subsection{Case when $n-1$ is not a power of $2$}

Here we take the surface $S=S(n-2,1,1,2,0,1)$. This surface is 
$\hp{1}{\cdot}\Lsurf(n-1,2)$.

\begin{wrapfigure}{r}{0pt}
  \includegraphics{anf.eps}
\end{wrapfigure}

The cusp width of $S$ is $2$, and $S$ has one vertical cylinder, hence
vertical cusp width $n$.  So $\Gamma$ contains $\hp{2}$ and $\vp{n}$.

Here $d = \frac{3}{16} (n-1) \prod_{p\mid n} p^{2\nu-2} (p^2-1)$.

The choice of $m = \ell/\!/2n = \ell/(2^\lambda \prod_{p\mid
n}p^\lambda)$ leads to $\delta = 2^{3\lambda-2} {\cdot} 3 {\cdot}
\prod_{p\mid n} p^{3\lambda-2} (p^2-1)$.

It follows that $d$ divides $\delta$ if and only if $(n-1)$ divides
$2^{3\lambda+2} {\cdot} \prod_{p\mid n} p^{3\lambda-2\nu}$.

Since $n$ is not some $2^{k} + 1$, $n - 1$ has odd prime factors; these
do not divide $n$, so $d$ does not divide $\delta$ and $\Gamma$ is not a
congruence subgroup.

\section{Noncongruence of $\Gamma_{B_n}$ for odd $n\geqslant 5$}

\subsection{A proof for most cases}
\label{sec:b:most}
Consider the one-cylinder surface $S$ having saddle connections of 
lengths $1$, $n-3$, $2$ on the top and $2$, $n-3$, $1$ on the bottom.

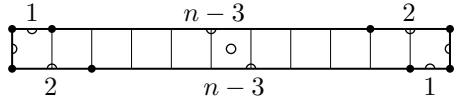
\begin{center}
\psfrag{1}{\Small $1$}
\psfrag{n}{\Small $n-3$}
\psfrag{2}{\Small $2$}
\includegraphics{b.eps}
\end{center}

The stabiliser of this surface contains $\hp{n}$ and $\vp{n}$.

Here $d = \frac{3}{16}(n-3)\prod_{p\mid n}p^{2\nu-2}(p^2-1)$.

The choice of $m = \ell /\!/ n$ leads to $\delta = 
\prod_{p\mid n}p^{3\lambda-2}(p^2-1)$.

Thus $d$ divides $\delta$ if and only if $3 (n-3)$ divides 
$16 \prod_{p\mid n} p^{3\lambda-2\nu}$.

Call an odd $n\geqslant 5$ ``bad'' if $3 (n-3)$ 
divides $16 \prod_{p\mid n} p^{3\lambda-2\nu}$.

As we are about to see, this is very rare, so that for ``most'' odd 
$n\geqslant 5$, $d$ does not divide $\delta$.

\subsection{The bad case}

If $n$ is such that $3 (n-3)$ divides $16\prod_{p\mid n} p^{3\lambda-2\nu}$,\\
\textbullet\ $n - 3$ is not a multiple of $2^5$;\\ 
\textbullet\ $n$ is a multiple of $3$ (and hence $(n - 3) \wedge n = 
3$);\\
\textbullet\ all odd prime factors of $n - 3$ divide $n$.

Combining these three remarks, we see the bad case is when $n - 3$ is
of the form $2^r {\cdot} 3^s$ with $1\leqslant r \leqslant 4$ and $1
\leqslant s$.

Thus the bad case consists of the four sequences $n_{r,s} =
2^r{\cdot}3^s+3$ for $r = 1$ to $4$ and $s\geqslant 1$, which have
exponential growth, hence zero density.

In particular, the discussion in \S\,\ref{sec:b:most}
proves the noncongruence of
$\Gamma_{B_n}$ when $n$ is out of these four sequences.

\subsection{First bad cases}
\label{sec:first:bad:cases}
Here we examine the first element of each of the four sequences, 
i.e.\ $n \in \{9,15,27,51\}$. We include the second element 
of the first sequence, i.e.\ $n = 21$.

\begin{wrapfigure}{r}{0pt}
\psfrag{1}{\Small $1$}
\psfrag{5}{\Small $5$}
\psfrag{n}{\Small $n-4$}
\includegraphics{b15.eps}
\end{wrapfigure}

Take $S = \Lsurf(5,n-4)$.  
Its horizontal cusp width is $n - 4$ and its vertical cusp width is 
$5$, so its stabiliser $\Gamma$ contains
$\hp{n-4}$ and $\vp{5}$.


Here $d = \frac{3}{16}(n-3)n^2\prod_{p\mid n}(1-\frac{1}{p^2})$.

Choosing $m = \ell /\!/ 5(n-4)$ leads to 
$\delta = \prod_{p\mid 5(n-4)}p^{3\lambda-2}(p^2-1)$.

The values of $d$ and $\delta$ for $n \in \{9,15,21,27,51\}$ are:

\medskip

\noindent
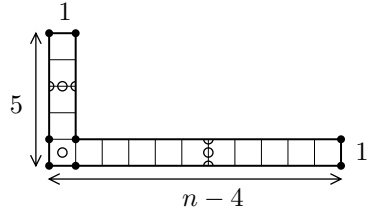
\centerline{%
\begin{tabular}{c@{\quad}c@{\quad}c@{\quad}c@{\quad}c@{\quad}c}
\hline
$n$ & $9$ & $15$ & $21$ & $27$ & $51$\\
\hline
$d$ & $3^4$ & $2^4 {\cdot} 3^3$ & $2^4 {\cdot} 3^4$ & $2^2 {\cdot}
3^6$ & $2^8 {\cdot} 3^4$\\
$\delta$ 
& $2^3 {\cdot} 3 {\cdot} 5$ 
& $2^6 {\cdot} 3^2 {\cdot} 5^2 {\cdot} 11$ 
& $2^8{\cdot}3^3{\cdot}5{\cdot}17$
& $2^7 {\cdot} 3^2 {\cdot} 5^4 {\cdot} 11 {\cdot} 23$ 
& $2^8 {\cdot} 3^2 {\cdot} 5^4 {\cdot} 23 {\cdot} 47$\\
\hline
\end{tabular}%
}

\medskip

In each case, we see by observing the power of $3$ in $d$ and $\delta$ 
that $d$ does not divide $\delta$.

\subsection{Remaining bad cases}

Here we will consider two surfaces $S_1$ and $S_2$ in orbit $B_n$, and
for each $S_i$ find some $k_i$ and $k'_i$ such that $\hp{k_i}$ and
$\vp{k'_i}$ are in the stabiliser $\Gamma_i$ of $S_i$ ($i\in\{1,2\}$).
The groups $\Gamma_1$ and $\Gamma_2$, being conjugate, have the same
index $d$ in $\Gamma(1)$ and the same level $\ell$.  Using $m_i = \ell
/\!/k_ik'_i$ will yield a $\delta_i$ for each $i\in\{1,2\}$ and we
will show that $d$ cannot divide both $\delta_1$ and $\delta_2$, 
implying that $\Gamma_{B_n}$ is not a congruence subgroup.

\begin{wrapfigure}{r}{0pt}
\psfrag{1}{\Small $1$}
\psfrag{2r}{\Small $2^r$}
\psfrag{3s}{\Small $3^s$}
\psfrag{3}{\Small $3$}
\includegraphics{bb1.eps}
\end{wrapfigure}

Take $S_1=S(1,2^r,3,3^s,1,0)$.

Its stabiliser contains $\hp{k_1}$ and $\vp{k'_1}$, with $k_1 =
3^s$ and $k'_1 = 2^r{\cdot}(2^r{\cdot}3+3)$.  Note that here $k'_1$ is
not the exact vertical cusp width of $S_1$, but a multiple of it.  
For $r = 1$, $2$, $3$, $4$, the value of $k'_1$ is respectively
$2{\cdot}3^2$, $2^2{\cdot}3{\cdot}5$, $2^3{\cdot}3^3$,
$2^4{\cdot}3{\cdot}17$.

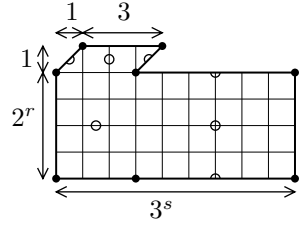
\begin{wrapfigure}{r}{0pt}
\psfrag{1}{\Small $1$}
\psfrag{2r}{\Small $2^r$}
\psfrag{3s}{\Small $3^s$}
\psfrag{3}{\Small $3$}
\includegraphics{bb2.eps}
\end{wrapfigure}

Take $S_2=S(3,2^r,1,3^s,0,0)$.


Its stabiliser contains $\hp{k_2}$ and $\vp{k'_2}$, with $k_2
= 3^s$ and $k'_2 = 2^r{\cdot}(2^r + 3)$; again $k'_2$ is not the exact
vertical cusp width, but a multiple of it.  It equals $2{\cdot}5$,
$2^2{\cdot}7$, $2^3{\cdot}11$, $2^4{\cdot}19$, respectively for $r =
1$, $2$, $3$, $4$.

Recall that $n = 2^r{\cdot}3^s + 3$, with $s \geqslant 2$.


Here, $d = \frac{3}{16}(n-3)n^2\prod_{p\mid n}(1-\frac{1}{p^2})$.
Since $3^2$ divides $(n-3)$, it does not divide $n$. Hence we can rewrite
$d = 2^{r-1} {\cdot} 3^{s+1+2\nu-2} \prod_{p\mid
\frac{n}{3}} p^{2\nu-2} (p^2-1)$.

%

The choice of $m_i = \ell /\!/ k_ik'_i$ leads to
$\delta_i = \prod_{p\mid k_ik'_i}p^{3\lambda-2}(p^2-1)$.


Given the values of $p^2 - 1$ for $p\in \{2, 3, 5, 7, 11, 17, 19\}$ 
(cf.\ table),


\medskip

\noindent
\centerline{%
\begin{tabular}{c@{\quad}c@{\quad}c@{\quad}c@{\quad}c@{\quad}c@{\quad}c@{\quad}c}
\hline
$p$ & $2$ & $3$ & $5$ & $7$ & $11$ & $17$ & $19$
\\
\hline
$p^2 - 1$ & $3$ & $2^3$ & $2^3 {\cdot} 3$ & $2^4 {\cdot} 3$ & $2^3 {\cdot}
3 {\cdot} 5$ & $2^5 {\cdot} 3^2$ & $2^3 {\cdot} 3^2 {\cdot} 5$\\
\hline
\end{tabular}%
}

\medskip

\noindent
the prime factors of $\delta_1$ and $\delta_2$ for each 
$r\in\{1,2,3,4\}$ are:

\medskip

\noindent
\centerline{%
\begin{tabular}{c@{\quad}c@{\quad}c@{\quad}c@{\quad}c}
\hline
$r$ & $1$ & $2$ & $3$ & $4$\\
\hline
$\delta_1$ & $2,3$ & $2,3,5$ & $2,3$ & $2,3,17$\\
$\delta_2$ & $2,3,5$ & $2,3,7$ & $2,3,5,11$ & $2,3,5,19$\\
\hline
\end{tabular}%
}

\medskip

%
%

If $d$ divides $\delta_1$ and $\delta_2$, we deduce that $\prod_{p\mid
\frac{n}{3}} p^{2\nu-2} (p^2-1)$ can have only $2$ and $3$ as prime
factors.  If this is the case, then $n$ has no square factor, and,
by Lemma \ref{lem:2:3:primes} (postponed to the end of the section), its prime factors are
in $\{3,5,7,17\}$.  The integers of the form
$3{\cdot}5^a{\cdot}7^b{\cdot}17^c$ with $a,b,c\in\{0,1\}$ are $3$,
$15$, $21$, $51$, $105$, $255$, $357$, $885$.  The only bad ones are
$15$, $21$, and $51$, and these were dealt with in
\S\,\ref{sec:first:bad:cases}.

To complete the proof of Theorem \ref{thm:noncongruence}, there 
remains only to prove:

\begin{Lemma}
\label{lem:2:3:primes}
If $p$ is prime and $p^2 - 1$ has no other prime factors than $2$ and 
$3$, then $p \in \{2,3,5,7,17\}$.
\end{Lemma}

This follows from the fact that $8$ and $9$ are the only two
consecutive nontrivial powers, a famous long-standing conjecture that
was recently proved by Mih\u{a}ilescu \cite{Mi}.

\newpage
\begin{ExternalTheorem}[Catalan's Conjecture]
The equation 

\noindent
\centerline{$x^u - y^v = 1,\ 
x>0,\ y>0,\ u>1,\ v>1$}

\noindent
has no other integer solution than $x^u = 3^2,\ y^v = 2^3$.
\end{ExternalTheorem}

\begin{proof}[Proof of the lemma]
%
%
By Catalan's Conjecture, consecutive powers of $2$ 
and $3$ are: ($1$, $2$); ($2$, $3$); ($3$, $4$); ($8$, $9$).
Suppose $(p-1)(p+1)$ has no other prime factors than $2$ and 
$3$. If $p$ is odd, then exactly one of $p-1$, $p+1$ is a multiple 
of $4$, and the other one is $2{\cdot}3^\alpha$.
If $\frac{p-1}{2} = 3^\alpha$, then either $\alpha = 0$, and $p = 3$, 
or $\alpha = 1$, and $p = 7$.
If $\frac{p+1}{2} = 3^\alpha$, then either $\alpha = 1$, and $p = 5$, 
or $\alpha = 2$, and $p = 17$.
\end{proof}

%
\flushleft

%

\end{document}